\documentclass{amsart}

\usepackage{amscd}
\usepackage{amsfonts}
\usepackage{amsmath}
\usepackage{amssymb}
\usepackage{amsthm}
\usepackage{amsxtra}
\usepackage{array}
\usepackage{fancybox}
\usepackage{fancyhdr}
\usepackage{graphics}
\usepackage{latexsym}
\usepackage{mathrsfs}
\usepackage{upref}
\usepackage{verbatim}
\usepackage[arrow, matrix, curve]{xy}

\newtheorem{lemma}{Lemma}
\newtheorem{theorem}[lemma]{Theorem}
\newtheorem{proposition}[lemma]{Proposition}

\newtheorem{question}[lemma]{Question}

\makeindex

\begin{document}
\title{A Property Of Local Cohomology Modules Of Polynomial Rings}
\author{Yi Zhang}
\address{Dept. of Mathematics, University of Minnesota, Minneapolis,
MN 55455}
\email{zhang397@umn.edu}
\thanks{NSF support through grant DMS-0701127 is gratefully acknowledged.}
\date{}

\begin{abstract}
Let $R=k[x_1,\cdots, x_n]$ be a polynomial ring over a field $k$ of
characteristic $p>0,$ and let $I=(f_1,\cdots,f_s)$ be an ideal of
$R.$ We prove that every associated prime $P$ of $H^i_I(R)$
satisfies $\text{dim}R/P\geqslant n-\sum\text{deg}f_i.$ In
characteristic 0 the question is open.
\end{abstract}
\maketitle

Throughout this paper, $R= k[x_1,\cdots, x_n]$ is the ring of
polynomials in $n$ variables over a field $k$ of characteristic
$p>0.$ The main result of this paper is the following:

\begin{theorem}\label{TAP}
Assume $I=(f_1,\cdots,f_s)$ is an ideal of $R$ such that
$\sum\text{deg}f_i<n.$ If $P$ is an associated prime of $H^i_I(R),$
then $\text{dim}R/P\geqslant n-\sum\text{deg}f_i.$
\end{theorem}

Whether the same result holds in characteristic 0 is still an open
question.

Before we give the proof of Theorem~\ref{TAP}, we fix some
notations. Denote the multi-index $(i_1,\cdots,i_n)$ by $\bar{i},$
especially $\overline{p^l-1}=(p^l-1,\cdots,p^l-1)$ where $l$ is a
positive integer. Since all the results in this paper concern the
vanishing of local cohomology modules and since extending the field
is faithfully flat, in the sequel, we can always enlarge the field
to make it perfect and infinite. Then $R$ is a free $R^{p^l}$-module
on the $p^{ln}$ monomials $e_{\bar{i}}=x^{i_1}_1\cdots x^{i_n}_n$
where $0\leqslant i_j<p^l$ for every $j.$ Let
$F:R\xrightarrow{r\mapsto r^{p^l}}R$ be the Frobenius homomorphism
and denote the source and target of $F$ by $R_s$ and $R_t$
respectively, that is $F:R_s\rightarrow R_t.$ There are two
associated functors
$$F^*:R_s\text{-mod}\rightarrow R_t\text{-mod}$$ such that
$F^*(-)=R_t\otimes_{R_s}-,$ and
$$F_*:R_t\text{-mod}\rightarrow R_s\text{-mod}$$ which is the
restriction of scalars. For each $R_s$-module $N,$ we have
$$F^*(N)=R_t\otimes_{R_s}N=\bigoplus_{\bar{i}}(e_{\bar{i}}\otimes_{R_s}N).$$
For each $f\in \text{Hom}_{R_t}(M,F^*(N)),$ define
$f_{\bar{i}}=p_{\bar{i}}\circ f:F_*(M)\rightarrow N,$ where
$$F^*(N)\xrightarrow{y\mapsto e_{\bar{i}}\otimes p_{\bar{i}}(y)}
e_{\bar{i}}\otimes_{R_s}N$$ is the natural projection to the
$\bar{i}$-component. There is a duality theorem in~\cite{gL09}:

\begin{theorem}\label{TD} (Theorem 3.3 in~\cite{gL09})
For every $R_t$-module $M$ and every $R_s$-module $N$ there is an
$R_t$-linear isomorphism
\begin{align*}
\text{Hom}_{R_s}(F_*(M),N) &\cong \text{Hom}_{R_t}(M,
F^*(N))\\
g_{\overline{p^l-1}}(-)&\leftarrow
(g=\oplus_{\bar{i}}(e_{\bar{i}}\otimes_{R_s}g_{\bar{i}}(-)))\\
g&\mapsto\oplus_{\bar{i}}(e_{\bar{i}}\otimes_{R_s}g(e_{\overline{p^l-1}-\bar{i}}(-))).
\end{align*}
\end{theorem}

\begin{proposition}\label{van} Assume $I=(f_1,\cdots,f_s)$ is an ideal of $R$ such that
$\sum\text{deg}f_i<n.$ Then $H^0_m(H^i_I(R))=0$ for every maximal
ideal $m.$
\end{proposition}

\begin{proof} Let
$K^\cdot(\underline{f}^t,R)$ be the Koszul cocomplex of $R$ on
$f^t_1,\cdots,f^t_s,$ that is
$$0\rightarrow R\xrightarrow{d^0}\bigoplus_{1\leqslant
\alpha\leqslant s}R_\alpha \xrightarrow{d^1} \bigoplus_{1\leqslant
\alpha_1<\alpha_2\leqslant
s}R_{\alpha_1,\alpha_2}\xrightarrow{d^2}\cdots
\xrightarrow{d^{s-1}}R_{1,\cdots,s}\rightarrow 0$$ where each
$R_{\alpha_1,\cdots,\alpha_j}$ is just a copy of $R$ indexed by the
tuple $(\alpha_1,\cdots,\alpha_j)$ and the differentials
$$d^j:\bigoplus_{1\leqslant \alpha_1<\cdots<\alpha_j\leqslant
s}R_{\alpha_1,\cdots,\alpha_j}\rightarrow \bigoplus_{1\leqslant
\alpha_1<\cdots<\alpha_{j+1}\leqslant
s}R_{\alpha_1,\cdots,\alpha_{j+1}}$$ are given by
$$(d^j(r))_{\alpha_1,\cdots,\alpha_{j+1}}=\sum^{v=j+1}_{v=1}(-1)^vf^t_v
r_{\alpha_1,\cdots,\hat{\alpha}_v,\cdots,\alpha_{j+1}}.$$ Let $l>1$
be an integer and $F:R_s\xrightarrow{r\mapsto r^{p^l}}R_t$ be the
Frobenius homomorphism. Let
$$\phi:H^i(K^\cdot(\underline{f},R))\rightarrow
F^*(H^i(K^\cdot(\underline{f},R)))
=H^i(F^*(K^\cdot(\underline{f},R)))$$ be the map induced by the
chain map
$$\phi^\cdot:K^\cdot(\underline{f},R)\rightarrow
F^*(K^\cdot(\underline{f},R)) =K^\cdot(\underline{f}^{p^l},F^*(R))$$
that sends each $R_{\alpha_1,\cdots,\alpha_j}$ to
$F^*(R_{\alpha_1,\cdots,\alpha_j})\cong
R_{\alpha_1,\cdots,\alpha_j}$ via multiplication by
$f^{p^l-1}_{\alpha_1}\cdots f^{p^l-1}_{\alpha_j}.$ By Proposition
1.11 in~\cite{gL97}, $H^i_I(R)$ is the direct limit of
$$H^i(K^\cdot(\underline{f},R)) \xrightarrow{\phi}
F^*(H^i(K^\cdot(\underline{f},R)))\xrightarrow{F^*(\phi)}
(F^*)^{2}(H^i(K^\cdot(\underline{f},R)))\xrightarrow{(F^*)^{2}(\phi)}\cdots.$$

Since $H^i(K^\cdot(\underline{f},R))$ is a subquotient of
$\bigoplus_{1\leqslant \alpha_1<\cdots<\alpha_i\leqslant
s}R_{\alpha_1,\cdots,\alpha_i},$ assume that
$H^i(K^\cdot(\underline{f},R))\subseteq \bigoplus_{1\leqslant
\alpha_1<\cdots<\alpha_i\leqslant s}R_{\alpha_1,\cdots,\alpha_i}/Q$
for a submodule $Q.$ Notice that $\phi$ is the multiplication by
$f^{p^l-1}_{\alpha_1}\cdots f^{p^l-1}_{\alpha_i}$ in the
$(\alpha_1,\cdots,\alpha_i)$ component. For each
$((g_{\alpha_1,\cdots,\alpha_i})+Q)\in \bigoplus_{1\leqslant
\alpha_1<\cdots<\alpha_i\leqslant s}R_{\alpha_1,\cdots,\alpha_i}/Q,$
suppose $l>\sum\text{deg}g_{\alpha_1,\cdots,\alpha_i},$ then
\begin{align*}
&\text{deg}(g_{\alpha_1,\cdots,\alpha_i}\cdot
f^{p^l-1}_{\alpha_1}\cdots
f^{p^l-1}_{\alpha_i})\\
\leqslant&\text{deg}g_{\alpha_1,\cdots,\alpha_i}+
(p^l-1)\cdot\text{deg}(f_1\cdots f_s)\\
<& l+(p^l-1)\cdot(n-1) \leqslant
n\cdot(p^l-1)\\
=&\text{deg}x_1^{p^l-1}\cdots x_n^{p^l-1},
\end{align*}
where the second inequality holds since
$\text{deg}g_{\alpha_1,\cdots,\alpha_i}<l$ and $\text{deg}(f_1\cdots
f_s)<n.$ In other words, all $g_{\alpha_1,\cdots,\alpha_i}\cdot
f^{p^l-1}_{\alpha_1}\cdots f^{p^l-1}_{\alpha_i}$ have zero
$e_{\overline{p^l-1}}=x_1^{p^l-1}\cdots x_n^{p^l-1}$ components in
$F^*(R_{\alpha_1,\cdots,\alpha_i}).$ Since
$$F^*(\bigoplus_{1\leqslant \alpha_1<\cdots<\alpha_i\leqslant
s}R_{\alpha_1,\cdots,\alpha_i}/Q)\cong (\bigoplus_{1\leqslant
\alpha_1<\cdots<\alpha_i\leqslant
s}F^*(R_{\alpha_1,\cdots,\alpha_i}))/F^*(Q),$$
$\phi((g_{\alpha_1,\cdots,\alpha_i})+Q)=((g_{\alpha_1,\cdots,\alpha_i}\cdot
f^{p^l-1}_{\alpha_1}\cdots f^{p^l-1}_{\alpha_i})+F^*(Q))$ has zero
$e_{\overline{p^l-1}}$ component in
$F^*(H^i(K^\cdot(\underline{f},R))).$ Therefore
$$\phi_{\overline{p^l-1}}((g_{\alpha_1,\cdots,\alpha_i})+Q)=p_{\overline{p^l-1}}((g_{\alpha_1,\cdots,\alpha_i}\cdot
f^{p^l-1}_{\alpha_1}\cdots f^{p^l-1}_{\alpha_i})+F^*(Q))=0.$$ Since
the corresponding morphism
$$\psi:F_*(H^i(K^\cdot(\underline{f},R)))\rightarrow
H^i(K^\cdot(\underline{f},R))$$ under Theorem~\ref{TD} is exactly
$\phi_{\overline{p^l-1}},$ we see that
$\psi((g_{\alpha_1,\cdots,\alpha_i})+Q)=0.$

Since $H^0_m(H^i(K^\cdot(\underline{f},R)))$ is of finite length, we
can choose $l$ big enough such that $\psi$ sends a $k$-basis of
$H^0_m(H^i(K^\cdot(\underline{f},R)))$ to 0. Hence $\phi=0$ on
$H^0_m(H^i(K^\cdot(\underline{f},R)))$ by Theorem~\ref{TD}.

Since local cohomology commutes with direct limit, $H^0_m(H^i_I(R))$
has the ``generating morphism", see \cite[Definition 1.9]{gL97},
$$\phi: H^0_m(H^i(K^\cdot(\underline{f},R)))\rightarrow
F^*(H^0_m(H^i(K^\cdot(\underline{f},R)))),$$ which is zero as has
just been shown. Therefore, $H^0_m(H^i_I(R))=0$ by Proposition 2.3
in~\cite{gL97}.
\end{proof}

\begin{proof}[Proof of Theorem~\ref{TAP}]
Let $h=\text{height}P.$ Suppose on the contrary we have
$\text{dim}R/P< n-\sum\text{deg}f_i.$ By Noether's normalization
lemma (\cite[Exercise 5.16]{AM69}), there are elements $y_1,\cdots,
y_{n-h} \in R$ which are algebraically independent over $k$ and such
that $k[y_1,\cdots,y_{n-h}]\cap P=0$ and $R/P$ is integral over
$k[y_1,\cdots,y_{n-h}].$ Moreover, since $k$ is infinite, these
$y_1,\cdots, y_{n-h}$ can be chosen to be linear combinations of
$x_1,\cdots, x_n.$ Therefore, we can assume without loss of
generality that $x_1,\cdots, x_n$ are such that $k[x_1,\cdots,
x_{n-h}]\cap P=0$ and $P$ is a maximal ideal in
$R'=K[x_{n-h+1},\cdots,x_n],$ where $K$ is the fraction field of
$k[x_1,\cdots,x_{n-h}].$ By Proposition~\ref{van}, we get
$H^0_P(H^i_I(R'))=0,$ which contradicts the assumption that $P$ is
an associated prime of $H^i_I(R).$
\end{proof}

When $i=s,$ the local cohomology module $H^s_I(R)$ has a simple
``generating morphism" \cite[Definition 1.9]{gL97}, namely,
$R/I\xrightarrow{}R/I^{[p]}.$ If $H^0_m(R/I)=0$ in the first place,
then we get $H^0_m(H^s_I(R))=0$ immediately. Therefore,
Proposition~\ref{van} suggests the following:

\begin{question}\label{q}
Let $R=k[x_1,\cdots,x_n]$ be the ring of polynomials such that $k$
is any field, and let $m$ be any maximal ideal. Assume
$I=(f_1,\cdots,f_s)$ is an ideal of $R$ such that
$\text{deg}f_1+\cdots+\text{deg}f_s<n.$ Is $H^0_m(R/I)=0?$
\end{question}

\begin{proposition}
Assume the answer to Question \ref{q} is positive. Then
$\text{proj.dim.}R/I\leqslant \sum\text{deg}f_i.$
\end{proposition}
\begin{proof}
We claim that there is a linear combination of variables with
coefficients in $k,$ which is $R/I$-regular. Indeed, assume all
linear combinations of variables are zero-divisors of $R/I.$ Choose
distinct elements $c_1,\cdots,c_n$ of the field $k$ which are not
roots of unity, and set $y_t=c^t_1x_1+\cdots+c^t_nx_n.$ Then any $n$
of these are linearly independent over $k.$ Since each $y_t$ is a
zero-divisor of $R/I,$ it is contained in some associated prime of
$R/I.$ Since the associated prime ideals of $R/I$ are finitely many,
one of them must contain more than $n$ $y_t$'s, and hence contains
the maximal ideal $m=(x_1,\cdots,x_n).$ This contradicts the
assumption $H^0_m(R/I)=0.$

By Auslander-Buchsbaum theorem, it suffices to prove that
$$\text{depth}R/I\geqslant n-\sum\text{deg}f_i.$$ If
$\sum\text{deg}f_i=n-1,$ then $\text{depth}R/I\geqslant 1$ by the
claim above. If $\sum\text{deg}f_i<n-1,$ after a linear change of
variables, assume that $x_n$ is $R/I$-regular. Modulo $x_n,$ the new
ideal $\bar{I}=(I+(x_n))/(x_n)$ satisfies $H^0_{\bar{m}}
(\bar{R}/\bar{I})=0$ by a positive answer to Question \ref{q}. Since $\text{depth}R/I
=\text{depth}\bar{R}/\bar{I}+1,$ we are done by induction.
\end{proof}

M. Stillman asked the following question:
\begin{question}\label{SQ} (\cite[Problem 3.14]{PS09})
Fix a sequence of natural numbers $d_1,\cdots,d_s.$ Does there exist
a number $q,$ such that $\text{proj.dim.}R/I\leqslant q$ when $R$ is
a polynomial ring (over any field $k$) and $I$ is an ideal with
homogeneous generators of degrees $d_l,\cdots,d_s?$
\end{question}

Thus a positive answer to Question~\ref{q} would also give a
positive answer to M. Stillman's question and would even provide an
upper bound on $q.$

This result is from my thesis. I thank my advisor Professor Gennady
Lyubeznik for his guidance, support and patience. I thank Professor
Craig Huneke for pointing out M. Stillman's Question~\ref{SQ} and
its connection with the main theorem of this paper.

\end{document}